\documentclass[english,11pt]{article}
\usepackage[T1]{fontenc}
\usepackage[latin9]{inputenc}
\usepackage{amsmath}
\usepackage{amssymb}

\usepackage{mathrsfs}
\usepackage{hyperref}
\usepackage{xy}
\hypersetup{colorlinks=true,citecolor=blue}
\usepackage{babel}

\begin{document}

\title{The Petersson norm of the Jacobi theta function}

\author{David Hansen}

\maketitle
Let $\theta(z)=\sum_{n\in\mathbf{Z}}e^{2\pi in^{2}z}$ be the Jacobi
theta function. This is a modular form of weight $1/2$ for the group
$\Gamma_{0}(4)$, and is well-known to be square-integrable; in fact,
it's the first interesting non-cuspidal but square-integrable automorphic
form. In this note we compute the norm\[
\left\Vert \theta\right\Vert ^{2}:=\int_{\Gamma_{0}(4)\backslash\mathfrak{H}}y^{\frac{1}{2}}|\theta(z)|^{2}\frac{dxdy}{y^{2}}.\]

\textbf{Theorem. }\emph{The Petersson norm of $\theta$ is $\left\Vert \theta\right\Vert ^{2}=4\pi$.}

Rather surprisingly, I have never seen this number calculated anywhere,
and I have seen at least one prominent researcher introduce it as
a kind of {}``fundamental constant'' in a paper. The problem is
that the constant term of $\theta$ prevents one from immediately
realizing $\left\Vert \theta\right\Vert ^{2}$ as the residue of a
Rankin-Selberg style integral. We get around this by a little trick.

Fix an arbitrary odd prime $p$, and consider the integral\[
I_{p}(s)=\int_{[0,1]\times\mathbf{R}_{>0}}y^{s+\frac{1}{2}}\left(|\theta(z)|^{2}-|\theta(p^{2}z)|^{2}\right)\frac{dxdy}{y^{2}}.\]
This converges absolutely for $\mathrm{Re}s>1$ and is easily calculated
as\begin{eqnarray*}
I_{p}(s) & = & 2\int_{\mathbf{R}_{>0}}y^{s-1/2}\sum_{n\geq1,\, p\nmid n}e^{-4\pi n^{2}y}\frac{dy}{y}\\
 & = & 2\cdot(4\pi)^{1/2-s}\sum_{n\geq1,\, p\nmid n}n^{1-2s}\int_{\mathbf{R}_{>0}}y^{s-1/2}\frac{dy}{y}\\
 & = & 2\cdot(4\pi)^{1/2-s}\Gamma(s-\tfrac{1}{2})(1-p^{1-2s})\zeta(2s-1).\end{eqnarray*}
On the other hand, the function $y^{\frac{1}{2}}\left(|\theta(z)|^{2}-|\theta(p^{2}z)|^{2}\right)$
is invariant under the group $\Gamma_{0}(4p^{2})$, so folding up
gives\[
I_{p}(s)=\int_{\Gamma_{0}(4p^{2})\backslash\mathfrak{H}}E_{4p^{2}}(z,s)y^{\frac{1}{2}}\left(|\theta(z)|^{2}-|\theta(p^{2}z)|^{2}\right)d\mu(z),\]
where $E_{4p^{2}}(z,s)=\sum_{\gamma\in\Gamma_{\infty}\backslash\Gamma_{0}(4p^{2})}\mathrm{Im}(\gamma z)^{s}$
is the usual nonholomorphic Eisenstein series and $d\mu(z)=\frac{dxdy}{y^{2}}$.
This series has a simple pole at $s=1$ with residue $\frac{3}{\pi}\cdot[\Gamma_{0}(1):\Gamma_{0}(4p^{2})]^{-1}=\frac{1}{2p(p+1)\pi}$.
Hence taking residues gives\begin{eqnarray*}
\mathrm{res}_{s=1}I_{p}(s) & = & \frac{1}{2p(p+1)\pi}\int_{\Gamma_{0}(4p^{2})\backslash\mathfrak{H}}y^{\frac{1}{2}}\left(|\theta(z)|^{2}-|\theta(p^{2}z)|^{2}\right)d\mu(z)\\
 & = & \frac{1}{2\pi}\int_{\Gamma_{0}(4)\backslash\mathfrak{H}}y^{\frac{1}{2}}|\theta(z)|^{2}d\mu(z)-\frac{1}{2p(p+1)\pi}\int_{\Gamma_{0}(4p^{2})\backslash\mathfrak{H}}y^{\frac{1}{2}}|\theta(p^{2}z)|^{2}d\mu(z)\\
 & = & \frac{1}{2\pi}\left\Vert \theta\right\Vert ^{2}-\frac{1}{2p^{2}(p+1)\pi}\int_{\Gamma_{0}(4p^{2})\backslash\mathfrak{H}}y^{\frac{1}{2}}|\theta(z)|^{2}d\mu(z)\\
 & = & \frac{1}{2\pi}(1-p^{-1})\left\Vert \theta\right\Vert ^{2},\end{eqnarray*}
where the third line follows from changing variables in the second
integral via the involution $z\to\frac{-1}{4p^{2}z}$ and the transformation
law $\theta(\frac{-1}{4z})=\sqrt{\frac{2z}{i}}\theta(z)$. But our
first computation gives \[
\mathrm{res}_{s=1}I_{p}(s)=2(1-p^{-1}),\]
and $p$ was arbitrary, so $\left\Vert \theta\right\Vert ^{2}=4\pi$.
$\square$ \[
\]
\[
\]

\end{document}